\documentclass[11pt]{article}
\usepackage{latexsym}
\usepackage{amsmath}
\usepackage{amsfonts}
\usepackage{mathrsfs}
\usepackage{graphics}

\newtheorem{lem}{Lemma}[section]
\newtheorem{thm}{Theorem}[section]
\newtheorem{cor}{Corollary}[section]

\newtheorem{exa}{Example}[section]

\newtheorem{assum}{Assumption}[section]
\newcommand{\Rs}{\mathbb{R}}

\newcommand{\trace}{\mbox{trace} }
\newcommand{\diag}{\mbox{diag} }

\newcommand{\bz}{{\bf 0} }

\newcommand{\Sn}{{\cal S}_n }

\newcommand{\Sh}{{\cal S}_H }
\newcommand{\Sc}{{\cal S}_C }

\newcommand{\KK}{{\cal K} }

\newcommand{\TT}{{\cal T} }

\newcommand{\bpr}{{\bf Proof.} \hspace{1 em}}
\newcommand{\epr}{ \\ \hspace*{4.5in} $\Box$ }
\newcommand{\beq}{ \begin{equation} }
\newcommand{\eeq}{ \end{equation} }
\newcommand{\bt}{ \begin{tabular} }
\newcommand{\et}{ \end{tabular} }

\begin{document}

\bibliographystyle{plain}
\title{A Remark on the Manhattan Distance Matrix of a Rectangular Grid } 
\vspace{0.3in}           
        \author{A. Y. Alfakih 
  \thanks{E-mail: alfakih@uwindsor.ca}
  \thanks{Research supported by the Natural Sciences and
          Engineering Research Council of Canada.}   
\\
\\
          Department of Mathematics and Statistics \\
          University of Windsor \\
          Windsor, Ontario N9B 3P4 \\
          Canada. } 
          
\date{\today}  
\maketitle


\vspace{0.1in}

\begin{abstract}
Consider the Quadratic Assignment Problem (QAP): given two matrices $A$ and $D$,
minimize $\{$ trace ($AXDX^T$) : $X$ is a permutation matrix $\}$. 
New lower bounds were obtained recently (Mittelmann and Peng \cite{mp10}) for the QAP
where $D$ is either 
the Manhattan distance matrix of a rectangular grid 
or the Hamming distance matrix of a hypercube.  
In this note, we show that the results in \cite{mp10,pmx10} extend to the case where $D$ is a
spherical Euclidean distance matrix, which includes the Manhattan distance matrix and the Hamming distance
matrix as special cases.    
\end{abstract}

\section{Introduction}

The Quadratic Assignment Problem (QAP)( see, e.g., \cite{bcpp98} and references therein)
is the problem:
\[ \min_{X \in \Pi} \trace (AXDX^T),  \]
 where $A$ and $D$ are two given matrices of order $n$ and $\Pi$
is the set of $n \times n$ permutation matrices. 
New lower bounds were obtained recently (Mittelmann and Peng \cite{mp10}) for the QAP
where $D$ is either 
the Manhattan distance matrix of a rectangular grid 
or the Hamming distance matrix of a hypercube. Mittelmann and Peng exploit the fact that
for such matrices $D$ there exists a scalar $\lambda$ such that $\lambda E - D$ is positive
semidefinite, where $E$ is the matrix of all 1's.  
In this note, we show that the results in \cite{mp10,pmx10} extend to the case where $D$
is a spherical Euclidean distance matrix (EDM), which includes the  
Manhattan distance matrix and the Hamming distance matrix as special cases.  

\subsection{Notation}
The positive semidefiniteness of a symmetric real matrix $A$ is denoted by
$A \succeq \bz$. The set of symmetric real matrices of order $n$ is denoted by
$\Sn$. $I_n$ denotes the identity matrix of 
order $n$ and  $e_n$ denotes the vector of all 1's in $\Rs^n$. 
$E_n$ denotes the $n \times n$ matrix of all 1's, i.e., $E_n=e_n e_n^T$.
The subscript $n$ will be omitted if the dimension of $I_n$, $e_n$ and $E_n$ 
is clear from the context.  
$A^{\dagger}$ denotes the Moore-Penrose inverse of $A$.
For two matrices $A$ and $B$, $A \otimes B $ denotes the 
Kronecker product of $A$ and $B$.
For a real number $x$, $\lceil x \rceil$ denote the ceiling of $x$.
Finally, $\diag(A)$ denotes the vector consisting of  
the diagonal entries of a matrix $A$.  

\section{Preliminaries}

The following well-known \cite{ber09} lemmas will be needed in some of the proofs below. 

\begin{lem} \label{lem1} 
Let $A$ and $B$ be two $n \times n$ real matrices. 
\begin{enumerate}
\item $(A^T)^{\dagger}$= $(A^{\dagger})^T$ := $ A^{\dagger T}$.
\item $(AA^T)^{\dagger}$= $A^{\dagger T} A^{\dagger}$.
\item If $A^TB=0$ and $BA^T=0$, then 
            $(A + B)^{\dagger} = A^{\dagger} + B^{\dagger}$. 
\item If $A^TB=0$ and $BA^T=0$, then 
            rank $(A + B)$ = rank $A$ + rank $B$. 
\item If $A$ has full column rank then $A^{\dagger}= (A^TA)^{-1}A^T$. 
\item Let $C$ and $Q$ be two real matrices of orders $k \times k$ and $n \times k$ respectively,
where $Q^TQ=I_k$. Then  
$( Q C Q^T )^{\dagger} = Q C^{\dagger} Q^T $. \label{5} 
\end{enumerate}
\end{lem}

\begin{lem} \label{lem2}
Let $A$ and $B$ be two $n \times n$ real matrices. 
\begin{enumerate}
\item $(A \otimes B)^{\dagger} = A^{\dagger} \otimes B^{\dagger}$. 
\item Let $A \succeq \bz$ and $B \succeq \bz$. Then $(A \otimes B) \succeq \bz$. 
\item $ rank(A \otimes B)$ = rank $(A)$ rank $(B)$. 
\end{enumerate}
\end{lem}

\begin{lem}[Generalized Schur Complement] \label{gsc}
Given the real symmetric partitioned matrix
\[
M= \left[ \begin{array}{cc} A & B \\ B^T & C \end{array} \right].
\]
Then $M \succeq \bz$ if and only if $C \succeq \bz$, $A - B C^{\dagger}B^T \succeq \bz$, and the null space
of $C$ is a subset of the null space of $B$.
\end{lem}

\section{Euclidean Distance Matrices (EDMs)} 

An $n \times n$  matrix $D =(d_{ij})$ is called a {\em Euclidean 
distance matrix (EDM)} if there exist points 
$p^1,\ldots, p^n$ in some Euclidean space $\Rs^r$ such that 
\[  d_{ij}= |\!| p^i - p^j |\!|^2 \;\; \mbox{ for all } i,j=1,\ldots,n,  \] 
where $|\!| \;\, |\!|$ denotes the Euclidean norm.
Moreover, the dimension of the affine span of the points  
$p^1,\ldots, p^n$ is called the {\em embedding dimension} of $D$.  
Without loss of generality, we make the following assumption.
\begin{assum}
The origin coincides with the
centroid of the points $p^1,\ldots,p^n$. 
\end{assum}
It is well known \cite{sch35,yh38} that a symmetric matrix $D$ whose diagonal entries
are all 0's is an EDM if and only if $D$ is negative semidefinite
on the subspace 
\[M := \{ x \in \Rs^n : e^T x = 0 \}, \] 
where $e$ is the vector of all 1's in $\Rs^n$.
It easily follows that the orthogonal projection on $M$ 
is given by 
\beq \label{defJ}
J:= I - e e^T/n= I-E/n.
\eeq   
Let $\Sn$ denote the set of symmetric $n \times n$ real matrices and let
$\Sh$ and $\Sc$ be two subspaces of $\Sn$ such that:
\[ 
\Sh = \{ A \in \Sn : \diag(A)=\bz \} \mbox{ and } 
\Sc = \{ A \in \Sn : Ae=\bz \}.
\]  
Furthermore, let  \cite{cri88} $\TT: \Sh \rightarrow \Sc$
and $\KK : \Sc \rightarrow \Sh$ be the two linear maps defined by 

\begin{eqnarray}  
\TT(D) & := &  - \frac{1}{2} \; J D J,  \label{defT} \label{defT} \\ 
\KK(B) & := &  \diag(B) \; e^T + e \; (\diag(B))^T - 2 \, B. \label{defK}  
\end{eqnarray}

Then it immediately follows that $\TT$ and $\KK$ are
mutually inverse between the two subspaces $\Sh$ and $\Sc$;  
and that $D$ is an EDM of embedding dimension
$r$ if and only if the matrix  
$\TT(D)$ is positive semidefinite of rank $r$ \cite{cri88}. 
Moreover, it is not difficult to show (see, e.g., \cite{alf06b,gow85,thw96})
that if $D$ is an EDM of embedding
dimension $r$ then rank ($D$) is equal to either $r+1$ or $r+2$.

Note that if $D$ is an EDM of embedding dimension $r$ then $\TT(D)$ is the Gram matrix 
of the points $p^1,\ldots,p^n$, i.e., $\TT(D)=PP^T$, where  $P$ is $n \times r$ and
${p^i}^T$ is the $i$th row of $P$. i.e.,
\beq \label{defP}
P=\left[ \begin{array}{c} {p^1}^T \\ \vdots \\ {p^n}^T \end{array} \right].
\eeq
$P$ is called a
{\em configuration matrix } of $D$. Also, note that $P^T e=\bz$ which is consistent with
our assumption that the origin coincides with the centroid of the points
$p^1,\ldots,p^n$. 

\subsection{Spherical EDMs}

An EDM $D$ is called {\em spherical} if  
the points $p^1,p^2,\ldots, p^n$ that generate $D$ lie on a 
hypersphere. Otherwise, $D$ is called
{\em non-spherical}.  
The following characterization of spherical EDMs is well known.

\begin{thm}[\cite{neu81,gow85,thw96}]   \label{sc} 
Let $D \neq \bz$ be a given $n \times n$ EDM with embedding dimension 
$r$, and let $P$ be a configuration matrix of $D$ such that $P^Te=\bz$. 
Then the following statements are equivalent. 
\begin{enumerate}
\item $D$ is a spherical EDM.
\item Rank ($D$) = $r+1$.   
\item The matrix $\lambda \; e e^T - D $ is positive semidefinite
for some scalar $\lambda$. 
\item There exists a vector $a$ in $\Rs^r$ such that: 
\beq \label{defa1}   Pa = \frac{1}{2} J \diag \,( \TT (D)),   \eeq
where $J$ is as defined in (\ref{defJ}).  
\end{enumerate}
\end{thm}

Two remarks are in order here. First, 
spherical EDMs can also be characterized in terms of Gale transform
\cite{aw02} which lies outside the scope of this paper. Second,
if $D$ is an $n \times n$ spherical EDM then the points that generate $D$ lie
on a hypersphere of center $a$ and radius $\rho  = (a^Ta + e^TDe/2n^2)^{1/2}$,
where $a$ is given in (\ref{defa1}) (see \cite{thw96}). Furthermore,
as the next theorem shows, the minimum value 
of $\lambda$ such that  $\lambda E - D \succeq \bz$ is closely related to
$\rho$. 

\begin{thm}[\cite{neu81,thw96}] \label{srho}  
Let $D$ be a spherical EDM generated by points that lie on a
hypersphere of radius $\rho$. Then $\lambda^* = 2 \rho^2$ is the minimum
value of $\lambda$ such that $\lambda E - D \succeq \bz$.  
\end{thm}
The following lemma gives an explicit expression for the radius $\rho$ which is
more convenient for our purposes. 
 
\begin{lem} \label{rho}  
Let $D$ be a spherical EDM generated by points that lie on a
hypersphere of radius $\rho$. Then 
\beq \label{defrho}  
\rho^2 = \frac{e^T D e}{2n^2} + \frac{e^T D (\TT(D))^{\dagger} D e}{4n^2}. 
\eeq  
\end{lem}

\bpr Let $D$ be a spherical EDM and let $Q=[\frac{e}{\sqrt{n}} \; V]$
be an orthogonal matrix, i.e., $J = VV^T$. Then  
$\lambda E - D \succeq \bz$ if and only if  
$Q^T ( \lambda E - D) Q \succeq \bz$. But  
\[ Q^T ( \lambda E - D) Q = \left[ \begin{array}{ll}
          \lambda n - e^TDe/n & -e^T D V / \sqrt{n} \\
          - V^T D e /\sqrt{n} & - V^T D V  \end{array} \right]. 
\]
Clearly $-V^T D V$ = $2 V^T \TT(D) V \succeq \bz $. 
Thus it follows from Lemma \ref{gsc} that $\lambda E - D \succeq \bz$ if and only if  
$\lambda n - e^T D e/n - e^T D V (V^T \TT(D) V)^{\dagger} V^T D e/2n \geq 0$ and
the null space of $V^TDV$ is a subset of the null space $e^T DV $.  
But  it follows from Property (\ref{5}) of Lemma \ref{lem1} that 
$V (V^T \TT(D) V)^{\dagger} V^T$ = $(VV^T \TT(D) V V^T)^{\dagger}$ = $\TT(D)^{\dagger}$
since $VV^T=J$ and since $(\TT(D) )e = \bz$.
Moreover, since $D$ is a spherical EDM, it follows that the 
null space of $V^TDV$ is a subset of the null space $e^T DV $ \cite{aw02}.  
Thus, $\lambda E - D \succeq 0$ if and only if  
$\lambda n - e^T D e/n - e^T D  (\TT(D))^{\dagger} D e/2n \geq 0$ 
and the result follows by Theorem \ref{srho}. 
\epr

Note that $(\TT(D))^{\dagger}$=$(PP^T)^{\dagger}$ = $P^{\dagger T} P^{\dagger}$. Thus
it follows from (\ref{defrho}) that
$\rho^2=a^Ta +e^TDe/2n^2 = e^T De / 2n^2 + e^TD P^{\dagger T} P^{\dagger}D e/4n^2$. Therefore,
\beq \label{defa2}
a = \frac{1}{2n} P^{\dagger}D e = \frac{1}{2n} (P^TP)^{-1} P^T D e ,
\eeq 
where $P$ is a configuration matrix of $D$. 

An interesting class of spherical EDMs is that of regular EDMs.
A spherical EDM $D$ is said to be {\em regular} (also called EDM of strength 1 \cite{neu80,neu81})
if $D$ is generated by points that lie
on a hypersphere centered at the centroid of these points \cite{ht93}. i.e.,
\[
n a = P^T e = \bz,
\]
since we assume that the centroid of the points $p^1,\ldots,p^n$ is located at the
origin. Therefore, if $D$ is a regular EDM then $\diag (\TT(D)) = \rho^2 e$. Hence, it follows from
(\ref{defK}) that    
\[
De = 2 n \rho^2 e,
\]
where $\rho^2 = e^TDe/2n^2$.
Thus we have the following characterization of regular EDMs
\begin{thm}[\cite{ht93,neu80,neu81}]
Let $D$ be an EDM. Then $D$ is regular if and only if $e$ is an eigenvector of $D$. 
\end{thm}

\section{Main Results}

\begin{thm} \label{thmed} 
Let $D_1$ be an $m \times m$ EDM of embedding dimension $r_1$, and let
$D_2$ be an $n \times n$ EDM of embedding dimension $r_2$. 
Then $D= E_m \otimes D_2+ D_1 \otimes E_n$ is an EDM  
of embedding dimension $r_1+r_2$,
where $E_m$ is the $m \times m$ matrix of all 1's and 
$\otimes$ denotes the Kronecker product. 
\end{thm}

\bpr
$I_{nm}= I_m \otimes I_n$ and $E_{nm}= E_m \otimes E_n$. Thus
\begin{eqnarray*} 
\TT(E_m \otimes D_2) & = & -\frac{1}{2} (I_m \otimes I_n - \frac{1}{nm} E_m \otimes E_n) 
  (E_m \otimes D_2) (I_m \otimes I_n - \frac{1}{nm} E_m \otimes E_n) \\
       & = & -\frac{1}{2} E_m \otimes (I_n - \frac{1}{n} E_n) D_2  
  (I_n - \frac{1}{n} E_n)  \\ 
       & = & E_m \otimes \TT( D_2) \succeq \bz. 
\end{eqnarray*} 
Similarly $\TT(D_1 \otimes E_n)$ = $\TT(D_1) \otimes E_n \succeq \bz$. Therefore,  
       $\TT(D)= E_m \otimes \TT( D_2)+ \TT(D_1) \otimes E_n \succeq \bz$. 
Hence $D$ is EDM.

On the other hand, we have from Lemma \ref{lem1} that
rank ($\TT( D)$) = rank $(E_m \otimes \TT( D_2)+ \TT(D_1) \otimes E_n)$ 
= rank $(E_m \otimes \TT( D_2))$ + rank $(\TT(D_1) \otimes E_n)$ 
= rank ($\TT( D_2)$) + rank ($\TT(D_1)$) = $r_2 + r_1$. 
Thus the embedding dimension of $D $ is equal to $r_1+r_2$. 
\epr

\begin{thm} \label{thmr} 
Let $D_1$ be a spherical EDM of order $m$ generated by points that lie
on a hypersphere of radius $\rho_1$, and 
let $D_2$ be a spherical EDM of order $n$ generated by points that lie
on a hypersphere of radius $\rho_2$. Then 
$D=E_{m} \otimes D_2 + D_1 \otimes E_{n}$ is a spherical EDM    
generated by points that lie on a hypersphere of radius 
$\rho= (\rho_1^2 +\rho _2^2)^{1/2}$. 
\end{thm}

\bpr $\lambda E_{nm}- D$ = $ E_m \otimes (\lambda_2 E_n - D_2)$ +    
$ (\lambda_1 E_m - D_1) \otimes E_n$ where  
$\lambda$ = $\lambda_1$ + $\lambda_2$. 
Since $D_1$ and $D_2$ are spherical
EDMs it follows that $ \lambda_2 E_n -D_2 \succeq \bz$ for $\lambda_2 = 2\rho_2^2$,     
and $ \lambda_1 E_n -D_1 \succeq \bz$ for $\lambda_1 = 2\rho_1^2$. Therefore,
for $\lambda = 2 \rho_1^2 + 2 \rho_2^2$ we have      
$\lambda E_{nm}- D \succeq \bz $ and hence $D$ is a spherical EDM generated
by points that lie on a hypersphere of radius 
$\rho \leq (\rho_1^2+\rho_2^2)^{1/2}$.
Next we show that $\rho = (\rho_1^2+\rho_2^2)^{1/2}$.

By Lemma \ref{rho} we have 
\begin{eqnarray*}
\rho^2 & = & \frac{e^T D e}{2(nm)^2} + 
                      \frac{e^T D (\TT(D))^{\dagger} D e}{4(nm)^2} \\ 
 & = & \frac{e_n^T D_2 e_n}{ 2n^2} +  \frac{e_m^T D_1 e_m}{2m^2} \\ 
 & & + \{(m e_m^T \otimes e_n^T D_2 + e_m^T D_1 \otimes n e_n^T )   
   (E_m \otimes \TT(D_2) + \TT(D_1) \otimes E_n)^{\dagger} ... \\  
& &\;\;\;\;\;\;\;\;\;\;\;...(m e_m \otimes D_2 e_n + D_1 e_m \otimes n e_n )\}/4(nm)^2 \\ 
 & = & \frac{e_n^T D_2 e_n}{2n^2} +  \frac{e_m^T D_1 e_m}{2m^2} \\ 
 & & + \{(m e_m^T \otimes e_n^T D_2 + e_m^T D_1 \otimes n e_n^T )   
(\frac{E_m}{m^2} \otimes (\TT(D_2))^{\dagger} + (\TT(D_1))^{\dagger} \otimes \frac{E_n}{n^2}) ... \\  
 & & \;\;\;\;\;\;\;\;\;\;... (m e_m \otimes D_2 e_n + D_1 e_m \otimes n e_n )\}/4(nm)^2\\ 
 & = & \frac{e_n^T D_2 e_n}{2n^2} +  \frac{e_m^T D_1 e_m}{2m^2} +  
    \{m^2 e_n^T D_2 (\TT(D_2))^{\dagger} D_2 e_n \\   
 & & +  n^2 e_m^T D_1 (\TT(D_1))^{\dagger} D_1 e_m \}/(4(nm)^2 \\    
 & = & \rho_2^2 + \rho_1^2 
\end{eqnarray*}
\epr

\begin{exa} \label{ex1}
Let $G_n=(g_{ij})$ be the Manhattan distance matrix of a rectangular grid consisting of
1 row and $n$ columns. Then 
\[  g_{ij} = | i - j|.
\]
Let $q^1, \ldots,q^n$ be the points in $\Rs^{n-1}$ such that 
the first $i-1$ coordinates of $q^i$ are 1's and the remaining
$n-i$ coordinates are 0's.  
Then it immediately follows that $q^1,\ldots,q^n$
generate $G_n$ and $q^1,\ldots, q^n$ lie on a hypersphere of radius 
$\rho= \frac{1}{2}(n-1)^{1/2}$ and
centered at $b=(1/2,1/2,\ldots,1/2)^T$. 
Thus  $G_n$ is a spherical EDM of embedding dimension $n-1$. 
Note that in this example, we don't assume that 
the origin coincides with the centroid of the points $q^1,\ldots,q^n$ in order
to keep the expressions of the $q^i$'s and $b$ simple.
\end{exa}

Now consider a rectangular grid of $m$ rows and $n$ columns.
Let $\hat{d}_{ij,kl}$ be the Manhattan distance between the grid point 
of coordinates $(i,j)$ and the grid point of coordinates $(k,l)$. 
Then 
\[ \hat{d}_{ij,kl}=|i-k|+|j-l|.
\] 

Let $s=i+n(j-1)$ for $j=1,\ldots,m$ and $i=1,\ldots,n$; and  
let $t=k+n(l-1)$ for $l=1,\ldots,m$ and $k=1,\ldots,n$.
Then
\[
j=\lceil s/n \rceil, i = s - n (\lceil s/n \rceil -1) \mbox{ and }
l=\lceil t/n \rceil, k = t - n (\lceil t/n \rceil -1). 
\]
  
Define the $nm \times nm$ matrix 
$D =(d_{st})$ such that $d_{st}= \hat{d}_{ij,kl}$. Then it follows \cite{mp10} 
that 
\beq \label{defMd}
 D = E_m \otimes G_n + G_m \otimes E_n,   
\eeq   
where $G_n$ is as defined in Example \ref{ex1}.
Equation (\ref{defMd}) follows since $(E_m \otimes G_n)_{st}$ = $|i-k|$ where $i = s - n (\lceil s/n \rceil -1)$
and $k = t - n (\lceil t/n \rceil -1)$; and since 
$(G_m \otimes E_n)_{kl}$ = $|j - l|$ where $j=\lceil s/n \rceil$ and $l=\lceil t/n \rceil$.

Thus we have the following two corollaries. 

\begin{cor}
Let $D$ be the $mn \times mn$ Manhattan distance matrix of a rectangular
grid of $m$ rows and $n$ columns. Then $D$ is a spherical EDM    
of embedding dimension $n+m-2$. Furthermore, the points that generate $D$
lie on a hypersphere of radius $\rho= \frac{1}{2}(n+m-2)^{1/2}$.  
\end{cor} 

\bpr
This follows from Theorems \ref{thmed}  and \ref{thmr}. 
\epr

\begin{cor} [Mittlemann and Peng {\cite[Theorem 2.6]{mp10}}]  
Let $D$ be the $mn \times mn$ Manhattan distance matrix of a rectangular
grid of $m$ rows and $n$ columns. Then $(n+m-2)E_{mn}/2 - D \succeq 0$.  
\end{cor} 

\begin{exa} \label{ex2}
Consider the Hamming distance matrix of
the $r$-dimensional hypercube $Q_r$ whose vertices are the $2^r$ points
in $\Rs^r$ with coordinates equal to 1 or 0. 
Let $p^1, \ldots, p^{2^r}$ be the vertices of $Q_r$   
and let $D=(d_{ij})$ be the $2^r \times 2^r$ matrix where $d_{ij}$ is
the Hamming distance between $p^i$ and $p^j$; i.e.,
$d_{ij}= \sum_{k=1}^r |p^i_k-p^j_k|$. Then $d_{ij}$ is also equal to 
$\sum_{k=1}^r (p^i_k-p^j_k)^2$ =$|\!|p^i-p^j|\!|^2$. Therefore, $D$ is an
EDM. Furthermore, the points $p^1,\ldots,p^{2^r}$ lie on a 
hypersphere, centered at their centoid, of radius $\rho = \frac{1}{2} r^{1/2}$. 
Thus $D$ is a regular EDM of embedding dimension $r$. 
\end{exa}


\end{document}